\documentclass{ifacconf}
\usepackage{amsmath}
\usepackage{amssymb}
\usepackage{graphicx}      
\usepackage{natbib}        

\newcommand{\newsec}[1]{ {\bf \noindent #1. \quad \hspace{-0.5cm}}}

\newcommand{\Ln}{\mathrm{Ln}}

\newcommand{\grad}{\mathrm{grad}\,}
\newcommand{\divv}{\mathrm{div}\,}

\usepackage{amsmath,amssymb}

\newtheorem{prp}[thm]{Proposition}
\usepackage[all]{xy}
\usepackage{psfrag}

\begin{document}
\begin{frontmatter}

\title{Reaction-Diffusion Systems as Complex Networks} 

\author[a]{Marko Seslija}
\author[c]{Jacquelien M.A. Scherpen}
\author[b]{Arjan van der Schaft}
\address[a]{Departement Elektrotechniek, Katholieke Universiteit Leuven, Kasteelpark Arenberg 10,
B-3001 Leuven, Belgium, e-mail:~marko.seslija@esat.kuleuven.be}
\address[c]{Department of Discrete Technology and Production Automation, Faculty of Mathematics and Natural Sciences, University of Groningen, Nijenborgh 4, 9747 AG Groningen, The Netherlands, e-mail: J.M.A.Scherpen@rug.nl}
\address[b]{Johann Bernoulli Institute for Mathematics and Computer Science, University of Groningen, Nijenborgh 9, 9747 AG Groningen, The Netherlands, e-mail:~A.J.van.der.Schaft@rug.nl}

\vspace{0.2cm}
\begin{abstract}                
The spatially distributed reaction networks are indispensable for the understanding of many important phenomena concerning the development of organisms, coordinated cell behavior, and pattern formation. The purpose of this brief discussion paper is to point out some open problems in the theory of PDE and compartmental ODE models of balanced reaction-diffusion networks.
\end{abstract}

\vspace{0.5cm}
\end{frontmatter}
\vspace{0.5cm}
\section{Introduction}
Inspired by the recent advances in modeling and analysis of reaction networks, in \cite{SeslijaRDCom} we have provided a geometric formulation of the reversible reaction networks under the influence of diffusion. Exploiting the graph knowledge of the underlying reaction network, we have shown that the obtained reaction-diffusion system is a distributed-parameter port-Hamiltonian system on a compact spatial domain. 

\vspace{0.2cm}

Numerical methods are essential tools for the understanding of many important dynamical aspects of these complicated distributed port-Hamiltonian models. While there has been a number of computational techniques that proposed discretizations of reaction-diffusion equations, the geometric structures they model are often \emph{lost} in the process. In \cite{SeslijaRDCom}, we have offered a spatially consistent discretization of the PDE system and, in a systematic manner, recovered a compartmental ODE model on a simplicial triangulation of the spatial domain. Exploring the properties of the Laplacian of the complex network defined over a simplicial manifold, we have characterized the space of equilibrium points and provided a result that guarantees the spatiotemporal consensus of a large class of balanced reaction networks.

\vspace{0.2cm}

After a brief summary of the PDE and ODE models of reaction-diffusion networks, we shall formulate a few open problems pertaining to these systems.

\section{Reaction-Diffusion Networks}\label{Sec2}
The dynamics of a \textbf{\emph{balanced reaction network}} involving $m$ chemical species (metabolites) takes the form
\begin{equation}\label{masterequation}
\dot{x} = - Z B \mathcal{K}(x^*) B^\textsc{t} \mathrm{Exp} \left(Z^\textsc{t} \mathrm{Ln}\left(\frac{x}{x^*}\right)\right), 
\end{equation}
where $x\in \mathbb{R}_{+}^m$ represents the concentrations vector; $Z$ is an $m \times c$ {\it\textbf{complex stoichiometric matrix}}, whose $\rho$-th column captures the expression of the $\rho$-th complex in the $m$ chemical species; $B$ is a $c \times r$ {\it \textbf{incidence matrix}} capturing the topology of the complex graph; $\mathcal{K}(x^*) $ is a $r \times r$ positive diagonal matrix of balanced reaction constants given as $\mathcal{K}(x^*) := \mathrm{diag} \big( \kappa_1(x^*), \cdots, \kappa_r(x^*) \big)$; $x^* \in \mathbb{R}^m_+$ is a {\emph{\textbf{thermodynamic equilibrium}}}. There is a close relation between the matrix $Z$ and the standard stoichiometric matrix $S$, which is expressed as $
S = ZB$. For more details of balanced reaction networks see \cite{AJReaction}.

The form (\ref{masterequation}) is the starting point for the analysis of balanced chemical reaction networks in this paper. We shall assume the validity of the \emph{global persistency conjecture,} which states that for a positive initial condition $x_0\in \mathbb{R}_{+}^m$, the solution $x$ of (\ref{masterequation}) satisfies: $\mathrm{lim}\,\mathrm{inf}_{t\rightarrow \infty}x(t)>0$. The global persistency conjecture recently was proven for the single linkage class case in \cite{Anderson}, but for the system (\ref{masterequation}) remains an open problem.

\newsec{Stability of Balanced Reaction Networks}~It follows that once a thermodynamic equilibrium $x^*$ is given, the set of {\it all} thermodynamic equilibria is described by the following proposition.

\begin{prp}[\cite{AJReaction}]\label{thermoequilibrium1}
Let $x^* \in \mathbb{R}^m_+$ be a thermodynamic equilibrium, then the set of {\it all} thermodynamic equilibria is given by 
\begin{equation}\label{equilibria}
\mathcal{E} := \{ x^{**} \in \mathbb{R}^m_+ \mid S^\textsc{t} \Ln \left(x^{**}\right) = S^\textsc{t} \Ln \left(x^{*}\right) \}.
\end{equation}
\end{prp}

Making use of the formulation of the dynamics of balanced reaction networks in (\ref{masterequation}), in \cite{AJReaction} it was shown that all equilibria of a balanced reaction network are actually thermodynamic equilibria, and thus given by (\ref{equilibria}). 

The Gibb's free energy associated to the reaction system is given by
\begin{equation}\label{eq:FreeEnergyGibbs}
G(x) =x^\textsc{t} \mathrm{Ln}\left(\frac{x}{x^*}\right) +  \left(x^* - x \right)^\textsc{t} {\mathbf{1}}_m,
\end{equation}
where $x^*$ is an equilibrium of the reaction network and $\mathbf{{1}}_m$ denotes a vector of dimension $m$ with all ones.

Exploiting the properties of the balanced weighted Laplacian matrix $B \mathcal{K}(x^*) B^\textsc{t} $ and employing $G$ as a Lyapunov function, \cite{AJReaction} showed that all the thermodynamic equilibria are, in fact, asymptotically stable.

\newsec{PDE Model}
When a well-mixed hypothesis is not reasonable, a more appropriate model for the reaction network (\ref{masterequation}) is that of reaction-diffusion equations. To that end, let $M$ be a \emph{compact} $n$-dimensional smooth \emph{Riemannian} manifold with boundary $\partial M$, representing the spatial domain. The port-Hamiltonian reaction-diffusion system given in terms of the disagreement vector $\frac{x}{x^*}$ is given as
\begin{equation}\label{eq:pHRDvec-calRel3RdNew}
\begin{split}
  \frac{\partial x}{\partial t}&=\divv \!\!\left(R_d(x)\grad\!\!\left( \frac{x}{x^*}\right) \right)+f(x)\\
  e_b&=\mathrm{Ln}\left( \frac{x}{x^*} \right) |_{\partial M}\\
  f_b&=R_d(x)\grad\! \!\left(\frac{x}{x^*}\right)\cdot \nu|_{\partial M}\,,
\end{split}
\end{equation}
where $R_d$ is a positive-semidefinite diagonal diffusion matrix, the so-called energy-diffusion matrix, and $f$ is given by the right-hand side of (\ref{masterequation}). The system (\ref{eq:pHRDvec-calRel3RdNew}) is accompanied with the appropriate smooth initial condition. 

\section{Complex Networks}
Following the exposition of \cite{SeslijaAutomatica}, let $K$ be a homological simplicial complex obtained by triangulation of the manifold $M$. Assuming that $K$ is well-centered, its circumcentric dual is $\star K=\star_\mathrm{i} K\times \star_\mathrm{b} K$, where $\star_\mathrm{i} K$ is the interior dual and $\star_\mathrm{b} K$ is the boundary dual, as explained in \cite{SeslijaAutomatica,SeslijaJGP}.

\begin{figure}
\centering
\vspace{-0.0cm}
      \hspace{0cm}\includegraphics[width=8.0cm]{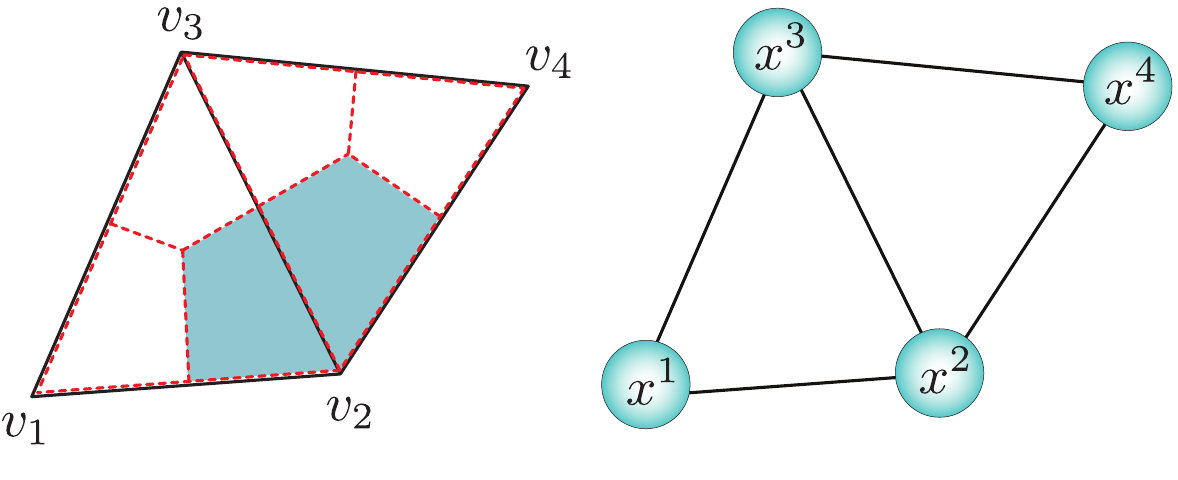}
      \vspace{-0.1cm}
  \caption{A simplicial complex $K$ consisting of two triangles. The dual edges introduced by the circumcentric subdivision are shown dotted. The state vector $x^j=\left(x_1^j,\ldots, x_m^j \right)^{\textsc{t}}$ is associated to the vertex $v_j$ for each $j\in\{1,\ldots, N\}$. The number of compartments for this example is $N=4$. The shaded region, the dual cell $\star_\mathrm{i} v_2$ of the vertex $v_2$, represents the compartment with the state $x^2$.}\label{fig:CH:RDcompart}
\end{figure}

The discrete analogue of an oriented manifold is an oriented simplicial complex, while differential forms are discretized as cochains. A $k$-cochain is a real-valued function on the $k$-simplices of $K$, which we will also call a \emph{discrete} $k$-form. Analogously, we define the space of discrete forms on $\star_\mathrm{i} K$ and $\star_\mathrm{b} K$. By $\Omega_d^k(K)$, $\Omega_d^k(\star_\mathrm{i} K)$, and $\Omega_d^k(\star_\mathrm{b} K)$ we denote the space of the primal $k$-cochains, the dual $k$-cochains, and the boundary dual $k$-cochains, respectively.

To each vertex of the primal mesh $K$ we associate reaction dynamics. That is, to a vertex $v_j$ we associate the state $x^j\in {\mathbb{R}}_+^m$. The geometric dual of $v_j$, $\star_\mathrm{i }v_j$, is the dual volume cell which represents the $j$-th compartment (see Figure~\ref{fig:CH:RDcompart}). The number of the compartments is $N=\mathrm{dim}\,\Omega_d^0(K)=\mathrm{dim}\, \Omega_d^n(\star_\mathrm{i} K)$.

By $X$ denote the concatenated vector
\begin{equation}
\label{eq:RDconcatenatedX}
X=\left( \left(x^1\right)^\textsc{t}, \dots, \left(x^N\right)^\textsc{t} \right)^\textsc{t},
\end{equation}
where $x^j\in {\mathbb{R}}_+^m$, and let
\begin{equation}
\label{eq:RDconcatenatedF}
F(X)=\left( f\left(x^1\right)^\textsc{t}, \dots, f\left(x^N\right)^\textsc{t} \right)^\textsc{t}
\end{equation}
be the vector field which describes the reaction dynamics of all compartments, with the reaction kinetics $f(x^j)=- Z B \mathcal{K}(x^*) B^\textsc{t} \mathrm{Exp} \left(Z^\textsc{t} \mathrm{Ln}\left(\frac{x^j}{x^*}\right)\right)$, $j=1,\dots, N$.

The \emph{\textbf{open compartmental model}} of the reaction-diffusion system (\ref{eq:pHRDvec-calRel3RdNew}) is given by
\begin{equation}\label{RD-DEC}
\begin{split}
\dot X&=-\left( \left(*_0\right)^{-1} \!\otimes I_m\right) \!\!\bigg( \Delta_d \frac{X}{X^*} -\left( \mathbf{tr}\otimes I_m\right)^\textsc{t}\hat f_b\bigg) +F(X)\\
e_b&= \left(\mathbf{tr}\otimes I_m\right) \frac{X}{X^*}\,,
\end{split}
\end{equation}
with the boundary flows $\hat f_b\in \left(\Omega_d^{n-1}(\star_\mathrm{b} K)\right)^{m}$ and the boundary efforts $e_b\in \left( \Omega_d^0(\partial K) \right)^m$. The symbol $\otimes$ represents the Kronecker product and $I_m$ is the identity matrix of dimension $m\times m$. The discrete Hodge operator $*_1:\Omega^1(K)\rightarrow \Omega^{n-1}( \star_\mathrm{i} K)$ is a diagonal matrix with the $k$-th entry being equal $|\star_\mathrm{i} \sigma_k^1|/|\sigma_k^1|$, where $\sigma_k^1$ is the primal edge with the dual $\star_\mathrm{i} \sigma_k^1$. The matrix $*_0$ is a diagonal matrix whose $k$-th element is $|\star_\mathrm{i} v_k|/|v_k|$. The Laplacian matrix of the simplicial complex is $\Delta_d=\left( \mathbf{d}\otimes I_m\right)^\textsc{t} \left(*_1 \otimes I_m\right) R_d(X) \left( \mathbf{d}\otimes I_m\right)$ with $\mathbf{d}$ being the discrete exterior derivative\footnote{The discrete exterior derivative $\mathbf{d}$ in this case is nothing but the transpose of the incidence matrix of the simplicial complex.}, $\mathbf{tr}$ is the trace operator\footnote{The trace operator $\mathbf{tr}$ is a matrix that isolates the members of a $0$-cochain vector assumed on the geometric boundary $\partial K$.}, $R_d(X)\geq \alpha I_{m N_e}$, $\alpha >0$, and $N_e$ is the number of edges of the primal mesh, i.e., $N_e=\mathrm{dim}\,\Omega_d^1(K)=\mathrm{dim}\, \Omega_d^{n-1}(\star_\mathrm{i} K)$. 

Furthermore, $\frac{X}{X^*}=\left( \left(\frac{x^1}{x^*}\right)^\textsc{t},  \left(\frac{x^2}{x^*}\right)^\textsc{t}, \dots,  \left(\frac{x^N}{x^*}\right)^\textsc{t}\right)^\textsc{t}$.

The total energy of the system, the sum of energies of all compartments, is
\[
G_d(X)=\sum_{j=1}^N G(x^j)V_{v_j}\,,
\]
where $G(x^j)$ is the free energy of the state $x^j$ and $V_{v_j}$ is the $n$-dimensional support volume obtained by taking the convex hull of the simplex $v_j$ and 
and its dual cell $\star_\mathrm{i} v_j$. Since $V_{v_j}=|v_j| |\star_\mathrm{i}v_j|= |\star_\mathrm{i}v_j|$, $j=1,\dots,N$, the total energy is $G_d(X)\!=\!\sum_{j=1}^{N} G(x^j) |\star_\mathrm{i}v_j|$.

\vspace{0.5cm}
\newsec{Compartmental Model}~Imposing the zero-flux boundary conditions, $\hat f_b=0$, leads to the closed compartmental model
\begin{equation}\label{eq:ClosedCompartmentalRD}
\begin{split}
\!\!\!\!\dot X&\!=\!-\left( \left(*_0\right)^{-1} \otimes I_m\!\right)  \Delta_d \frac{X}{X^*}  +F(X),\!\!
\end{split}
\end{equation}
with a positive initial condition $X(0)=X_0\in \mathbb{R}_+^{mN}$.

In \cite{SeslijaRDCom}, we have shown that for the system (\ref{eq:ClosedCompartmentalRD}) the positive orthant $\mathbb{R}_+^{mN}$ is forward invariant. In order to exclude the existence of possible boundary equilibria, given $X_0\in \mathbb{R}_+^{mN}$, we assume that all the trajectories $t\mapsto X(t)$ of (\ref{eq:ClosedCompartmentalRD}) satisfy: $\mathrm{lim}\,\mathrm{inf}_{t\rightarrow \infty}X(t)>0$.

In the absence of the diffusion terms, the dynamics of the spatially discrete systems in (\ref{eq:ClosedCompartmentalRD}) are decoupled, and as such coincide with the dynamics of the balanced reaction system (\ref{masterequation}). In this scenario all the compartments exhibit asymptotically stable dynamics, but the steady states of all the compartments, in general, are not identical. The following theorem shows that the compartmental model (\ref{eq:ClosedCompartmentalRD}) is asymptotically stable with the spatially uniform steady state. 

\begin{thm}[\cite{SeslijaRDCom}]\label{tm:LyapRDcompartment}
Consider the compartmental model of balanced mass action reaction network given by (\ref{eq:ClosedCompartmentalRD}). For every initial condition $X(0)\in \mathbb{R}_+^{mN}$, the species concentrations $x^1, \dots, x^N$ as $t \to \infty$ converge to $x^1=\cdots=x^N\in\mathcal{E}$.
\end{thm}

\section{Open Problems}
\newsec{Global Persistency Conjecture}
Proving the global persistency conjecture for balanced reaction networks remains an important unsolved problem with many significant mathematical and biochemical consequences.

\newsec{Stability of PDE reaction-diffusion networks}
The existence of solutions for the system (\ref{eq:pHRDvec-calRel3RdNew}) is a complex issue. The papers \cite{Morgan,Fitzgibbon} do provide a working framework for the systems with separable Lyapunov functions. Furthermore, according to \cite{Fitzgibbon}, the system \emph{does not }generate spatial patterns. However, proving the spatial uniformity of asymptotic behavior of the balanced reaction-diffusion systems involves Krasovskii-LaSalle-type of arguments. These arguments in turn require the precompactness and the global boundedness of classical solution. Showing the global boundedness of classical solution of the semilinear balanced eaction-diffusion networks in the presence of Neumann boundary conditions is a challenging issue.

\newsec{Pattern Formation}
It is well-known that adding diffusion to the reaction system can generate behaviors absent in the ODE case. This primarily pertains to the problem of diffusion-driven instability which constitutes the basis of Turing's mechanism for pattern formation \cite{Turing}, \cite{MurrayBio}. Here, the port-Hamiltonian perspective permits us to
draw immediately some conclusions regarding passivity of reaction-diffusion systems, but also to claim the spatial uniformity of the asymptotic behavior of the compartmental model.

Given a reaction-diffusion system, an important question is whether it is possible to generate spatial pattern by manipulating the boundary variables of the boundary control problem. Since both the smooth (\ref{eq:pHRDvec-calRel3RdNew}) and the discretized model (\ref{RD-DEC}) assume the port-Hamiltonian form, many elaborate schemes, ranging from passivity-based to optimal control, could potentially be applied for control of reaction-diffusion networks. Construction and study of such control strategies for reaction-diffusion networks can improve our quantitative understanding of pattern formation, but also may foster applications in bioengineering.

{

}


\begin{thebibliography}{99}
\providecommand{\natexlab}[1]{#1}
\providecommand{\url}[1]{\texttt{#1}}
\expandafter\ifx\csname urlstyle\endcsname\relax
  \providecommand{\doi}[1]{doi: #1}\else
  \providecommand{\doi}{doi: \begingroup \urlstyle{rm}\Url}\fi
\bibitem[Anderson, 2011]{Anderson}Anderson, D.F. (2011) ``A proof of the Global Attractor Conjecture in the single linkage class case,'' SIAM J. Appl. Math., Vol. 71, No. 4.

\bibitem[Fitzgibbon~{\em et al.}, 1997]{Fitzgibbon} Fitzgibbon, W.B., Hollis, S.L., Morgan, J.P. (1997)``Stability and
Lyapunov functions for reaction-diffusion systems,'' \emph{SIAM
Journal on Mathematical Analysis}, vol. 28, no. 3, pp. 595--610.

\bibitem[Morgan, 1991]{Morgan} Morgan, J. (1991) ``Global Existence for Semilinear Parabolic Systems on One-dimensional Bounded Domains,'' Rocky Mountain Journal of Mathematics, vol. 21, no. 2, 1991.

\bibitem[Murray, 2003]{MurrayBio} Murray, J. (2003) \emph{Mathematical Biology,} 3rd edition, Berlin: Springer-Verlag, 2003.

\bibitem[van der Schaft~{\em et al.}, 2012]{AJReaction} van der Schaft, A.J., Rao, S., Jayawardhana, B. (2012) ``On the Mathematical Structure of Balanced Chemical Reaction Networks Governed by Mass Action Kinetics,'' to appear in \emph{SIAM Journal on Applied Mathematics}, available at http://arxiv.org/abs/1110.6078v1.

\bibitem[Seslija~{\em et al.}, 2012a]{SeslijaJGP} Seslija, M., van der Schaft, A.J., Scherpen,  J.M.A.  (2012a) ``Discrete Exterior Geometry Approach to Structure-Preserving Discretization of Distributed-Parameter Port-Hamiltonian Systems," \emph{Journal of Geometry and Physics,} Volume 62, Issue 6, Pages 1509--153.

\bibitem[Seslija~{\em et al.}, 2012b]{SeslijaAutomatica} Seslija, M., Scherpen, J.M.A., van der Schaft, A.J. (2012b) ``Explicit Simplicial Discretization of Distributed-Parameter Port-Hamiltonian Systems,'' {arxiv.org/abs/1208.3549}, Submitted to \emph{Automatica} in August 2012.

\bibitem[Seslija~{\em et al.}, 2012c]{SeslijaRDCom} Seslija, M., van der Schaft, A.J., Scherpen, J.M.A. (2012c) ``Hamiltonian Perspective on Compartmental Reaction-Diffusion Networks,'' Submitted to \emph{Automatica} in December 2012.

\bibitem[Turing, 1952]{Turing} Turing, A.M. (1952) ``The chemical basis of morphogenesis,'' \emph{Philosophical trasactions of Royal Society of London}, Series B, Biological Sciences, Volume 237, Issue 641, pp. 37--72.

\end{thebibliography}
\end{document}